\documentclass[11pt, leqno]{amsart}
\usepackage{amsthm, amsmath}
\usepackage{amssymb} 
\usepackage{amscd}
\usepackage[curve,matrix,arrow,frame,tips]{xy}
\usepackage{verbatim}
\usepackage{graphicx}

\newtheorem{thm}{Theorem}[section]

\newtheorem{Remarksnumb}[thm]{Remarks}

\newtheorem{conjecture}[thm]{Conjecture}

\newcounter{ex}[section]

\newcommand{\cal}{\mathcal}

\newcommand{\E}{{\mathcal E}}

\newcommand{\G}{{\bf G}}

\newcommand{\C}{{\bf C}}

\newcommand{\Gg}{{\mathcal G}}

\newcommand{\Z}{{\bf Z}}
\newcommand{\F}{{\mathcal F}}
\newcommand{\ti}{\tilde}
\newcommand{\Spec}{{\rm Spec }\, }
 \renewcommand{\O}{{\mathcal O}}

\newcommand{\M}{{\mathcal M}}
\renewcommand{\L}{{\mathcal L}}

\newcommand{\lon}{{\longrightarrow}}

\newcommand{\Pg}{{\mathcal G}}

\newcommand{\loniso}{\buildrel \sim \over \lon}

\newcommand{\pione}{{\bf P}^1}

\def\thfill{\null\nobreak\hfill}

\def\endproof{\thfill\vbox{\hrule
  \hbox{\vrule\hbox to 5pt{\vbox to 5pt{\vfil}\hfil}\vrule}\hrule}}

\renewcommand{\P}{{\cal P}}

\numberwithin{equation}{section}

\begin{document}

\title[$\mathcal G$-bundles]{Some questions about $\mathcal G$-bundles on curves}
\author[G. Pappas]{G. Pappas}
\thanks{*Partially supported by NSF
grant DMS08-02686}
\address{Dept. of
Mathematics\\
Michigan State
University\\
E. Lansing\\
MI 48824-1027\\
USA} \email{pappas@math.msu.edu}
\author[M. Rapoport]{M. Rapoport}
\address{Mathematisches Institut der Universit\"at Bonn,
Beringstrasse 1\\ 53115 Bonn\\ Germany.}
\email{rapoport@math.uni-bonn.de}

\date{\today}

\maketitle

\bigskip
\bigskip

\section{ }

The purpose of this note is to discuss the geometry of  moduli 
stacks of various types of bundles over a  
curve. We suggest that 
the main elements of the theory of moduli of $G$-bundles for a constant  
  reductive algebraic group $G$ as developed by Beauville, Laszlo, Faltings and other
authors should extend to a theory of moduli of $\Pg$-torsors for a
large class of algebraic group schemes $\Pg$ 
which are not necessarily constant over the curve.
The class we consider is that of smooth group schemes over the curve 
with   reductive generic fiber which have the property that
each place of the curve the completion of the group scheme
is a ``parahoric group scheme" of the type constructed by Bruhat-Tits.
In addition to the classical case above, the corresponding moduli stacks
include the moduli of parabolic $G$-bundles but also 
other interesting examples such as the moduli of 
Prym line bundles (Prym varieties) or moduli of bundles together
with (not always  perfect)
symplectic, orthogonal or hermitian pairings.
Our approach uses the theory of loop groups.

In [PR], we introduced and studied the loop group attached to a linear algebraic group over a Laurent series field $k((t))$
where $k$ is an algebraically closed field. To a (connected) reductive
algebraic group $H$ over $k((t))$ there is associated the ind-group scheme $LH$ over $k$, with points with values in a $k$-algebra $R$ equal to $H(R((t)))$. If $P$ is a parahoric subgroup of $H(k((t)))$, Bruhat and Tits have associated to $P$ a smooth group scheme with connected fibers over $\Spec(k[[t]])$, with generic fiber $H$ and with group of $k[[t]]$-rational points equal to $P$. Denoting by the same symbol $P$ this group scheme, there is associated to it a group scheme $L^+P$ over $k$, with points with values in a $k$-algebra $R$ equal to $P(R[[t]])$. The $fpqc$-quotient $\F_P=LH/L^+P$ is representable by an ind-scheme, and is called the {\it partial affine flag variety} associated to $P$. In [PR] we studied these affine flag varieties
and obtained results about some of their basic structural properties. In particular, we  showed
\begin{itemize}
\item [1.)] $\pi_0(LH)=\pi_0(\F_P)=\pi_1(H)_I \ .$ Here $\pi_1(H)$ denotes the algebraic fundamental group of $H$ in the sense of Borovoi, and $I={\rm Gal}(\overline {k((t))}/k((t)))$ the inertia group. 
\item  [2.)] If $H$ is semi-simple and splits over a tamely ramified extension of $k((t))$, and if $({\rm char}(k), |\pi_1(H)|)=1$, then $LH$ and $\F_P$ are reduced ind-schemes. 
\end{itemize}
In the case when $H$ comes by extension of scalars from a constant algebraic group $H_0$ over $k$,   these properties and more have been shown in Faltings' paper \cite{FalJEMS}  (and much of it was known before, thanks to the work of Beauville, Laszlo, Sorger, Kumar, Littelmann, Mathieu, and others, comp. the references in \cite{FalJEMS} and  [PR]). In \cite{FalJEMS}, Faltings goes on to use these local results to prove global results on the moduli space of $H_0$-bundles on a smooth projective curve over $k$, in particular about its Picard group. The main tool
is the ``uniformization theorem" [BL], [DS], that expresses the
moduli stacks (for semisimple groups) as a certain quotient 
of the affine Grasmannian for $H$. In the present note we present some conjectures on how to generalize these results in the framework of [PR]. As it turns out, the results of Laszlo and Sorger in [LS] can be interpreted as a confirmation in special cases of our predictions. 

 After an older version of this paper was circulated,  
 Heinloth posted the preprint [He] where he proves a good part of these conjectures. 
 We hope that there is still some interest in our paper and that progress can be made in answering the rest of these
 questions. We also hope that our point of view can be a useful framework in generalizing 
 the enormous body of results for split groups to this more general case.

 For example, the correct extension of the Verlinde formula \cite{S1} in this context is still a mystery to us, although the articles \cite{BFS, FS, Sz} may contain ideas that could be
useful for this question. In any case, it would be interesting to
understand the relation between the spaces of ``conformal blocks" 
in our set-up of $\mathcal G$-bundles and the
spaces of conformal blocks for orbifold models of conformal
field theories that appear in these papers.

We thank G. Faltings, C. Sorger and J. Heinloth for helpful discussions. We also thank J. Heinloth and E. Looijenga for correcting some statements in the 
original version of this paper.

\section{ }

Let $k$ be an algebraically closed field, and let $X$ be a smooth connected projective curve over $k$. Let $\mathcal G$ be a smooth affine group scheme over $X$ with all fibers connected. In addition, we assume that the generic fiber $\mathcal G_{\eta}$ is a connected reductive group scheme over $K=k(X)$, and that for every $x\in X(k)$, denoting by $\mathcal O_x$ the completion of the local ring of $X$ at $x$ and by $K_x$ its fraction field, $\Gg(\O_x)\subset \Gg(K_x)$
is a parahoric  subgroup of $\mathcal G_{\eta_x}(K_x)$ in the sense of [BTII], see also [T]. 
We will call such a $\Gg$ a {\it parahoric group scheme over} $X$. Recall that by [BTII],
given a parahoric subgroup 
 $P_x\subset \mathcal G_{\eta_x}(K_x)$
 there is a unique affine smooth group scheme $\mathcal G_{P_x}$ over $\O_x$ with the following propreties: 
Its generic fiber is $\mathcal G_{\eta_x}$, it has connected special fiber and satisfies $\mathcal G_{P_x}(\O_x)=P_x$.

Let $\mathcal M_{\mathcal G/X}$ denote the stack of $\mathcal G$-torsors on $X$. 
The usual arguments show that this is a smooth (Artin) algebraic stack
over $k$. We are going to state four conjectures on the geometry of  $\mathcal M_{\mathcal G/X}$ but in this section we will first  discuss several examples.

\subsection{}\label{constant}
 
 Let $G$ be connected reductive group scheme over $k$. Then $G\times_{\Spec(k)}X$ is an example of the kind of group schemes we  consider. This is the case of a {\it constant group scheme}. 
 
 We may generalize this as follows. 
Let $x\in X(k)$. Then the parahoric subgroups in $G(K_x)$ contained in $G(\O_x)$ are in one-to-one correspondence with  the parabolic subgroups of $G$. More precisely, if $P\subset G$ is a parabolic subgroup, then the corresponding parahoric subgroup $\P$ is equipped with a morphism of group schemes over $\Spec \O_x$,
\begin{equation}
\P\rightarrow  G\times_{\Spec(k)}\Spec( \O_x)
\end{equation}
which in the generic fiber is the identity of $\mathcal G_{\eta_x}$ and which in the special fiber has image equal to $P$. 

Suppose now that $\mathcal G$ is a group scheme equipped with a morphism $\mathcal G\to G\times_kX$ which, when localized at $x$  is of the previous nature for all $x\in X(k)$. Hence there is a finite set of points $\{x_1,\ldots, x_n\}$ such that this morphism is an isomorphism outside this finite set, and  parabolic subgroups $P_1,\ldots, P_n$ such that the localization of $\mathcal G$ at $x_i$ corresponds to $P_i$ in the sense explained above.   Then there is an equivalence of categories between the category of $\mathcal G$-torsors on $X$ and the category of $G$-torsors on $X$ with quasi-parabolic structure of type $(P_1,\ldots, P_n)$ with respect to $(x_1,\ldots, x_n)$, in the sense of [LS].  For such group schemes 
some of the questions here have been considered in the literature, although not always in our formulation.  

\subsection{}\label{prym}

Let $S$ be a torus over $k(X)$; then the connected Neron model $\Gg={\mathcal S}^0$ of
$S$ over $X$ is another example of a parahoric group scheme. This kind of $\Gg$-bundle occurs in various other contexts that we mention here briefly.

Suppose that $\pi: Y\to X$ is an irreducible  finite   flat
and generically unramified covering.
Then $k(Y)/k(X)$ is a finite separable field extension
and we can take $S$ to be the torus ${\rm Res}_{k(Y)/k(X)}(\mathbb G_m)$, with parahoric extension $\Gg$ over $X$ equal to ${\rm Res}_{Y/X}(\mathbb G_m)$. Then a $\Gg$-bundle on $X$ is simply a line bundle $\L$ on $Y$. By associating to $\L$ its direct image $\pi_*(\L)$, we obtain a vector bundle of rank $n$  on $X$, where  $n=[k(Y):k(X)]$. This construction of vector bundles on $X$ is analyzed in [BNR]. If $Y$ is the  curve associated  in the sense of [BNR], \S 3 to a line bundle $\M$ and sections $\{ s_i\in \Gamma(X, \M^i)\mid i=1,\ldots, n\}$, then the vector bundles obtained by this construction have a canonical Higgs structure (wrt. $\M$), such that $Y$ is the associated spectral curve. 

We also mention the following variant, cf.\ \cite{Do}, \cite{DoGai}.   Suppose that $G$ is a connected reductive group with  maximal torus $T$, normalizer $N$ of $T$, and 
Weyl group $W$. Let   $\pi: Y\rightarrow X$ be a unramified Galois covering with Galois group $W$. Assume that the characteristic of $k$ does not divide the order of $W$. We can consider the 
group scheme
$$
{\mathcal G}=({\rm Res}_{Y/X}(T\times_k Y))^W
$$
 on $X$, where 
$W$ acts diagonally on $T\times_k Y$. Then, thanks to our assumption on the characteristic of $k$, $\Gg$ is a parahoric group scheme on $X$, cf. [E], Thm. 4.2. Each $\mathcal G$-torsor 
over $X$ gives an element in ${\rm H}^1(Y, T)^W={\rm Hom}_{W}(X^*(T), {\rm Pic}(Y))$
(here $W$ acts on both source and target). 
In general, if $\M$ is a $T$-torsor over $Y$ whose class belongs to
 ${\rm H}^1(Y, T)^W$, we can consider the group  $N_\M$ 
of automorphisms of $\M$ which commute with the action on $Y$ of some $w\in W$.
This affords an extension
$$
1\rightarrow T\rightarrow N_\M\rightarrow W\rightarrow 1\ .
$$
Suppose now that there is a $T$-torsor  $\M_0$ in $({\rm H}^1(Y, T))^W$ such that the corresponding extension $N_{\M_0}$ is isomorphic to
the extension given by the normalizer $N$ of $T$ in $G$. 
Then for each 
$\Gg$-torsor  on $X$, corresponding to the $T$-torsor $\L$ on $Y$,   $\pi_*(\L\otimes_{\O_Y}\M_0)$ gives an $N$-torsor over $X$
that can be induced
to give a $G$-bundle on $X$ (\cite{Do}, \cite{DoGai}). This $G$-bundle is an ``abstract" Higgs bundle with unramified {\sl cameral cover} $\pi: Y\to X$, loc. cit. This is a protypical result
in the theory of Higgs bundles and the Hitchin fibration. Here the (not  precise) catch-phrase is that the sufficiently generic fibers of the Hitchin map are --non-canonically-- isomorphic to   moduli varieties of $\Gg$-torsors for a suitable commutative $\Gg$ (a version of the above works even when $\pi$ is ramified, see \cite{Do}, \cite{DoGai}, \cite{Ngo}).

As an example consider the case $G={\rm SL}_2$. Then  $W=\Z/2\Z$, and $T=\mathbb G_m$ with $W$ acting 
by inversion.  Suppose that  ${\rm char}(k)\neq 2$ and that $\pi$ is an unramified double
cover with involution $\sigma$. The  parahoric group scheme $\Gg$ above is then the kernel 
 of the norm
$
{\rm Norm}_{Y/X}: {\rm Res}_{Y/X}(\mathbb G_m)\rightarrow \mathbb G_m
$.
In this case, the above amounts to a Prym construction which goes as follows (cf. \cite{Do} 5.2).
We can see that
 $\Gg$-torsors over $X$ 
are given by line bundles $\L$ on $Y$ such that ${\rm Norm}_{Y/X}(\L)$
is trivial.  
We can also see (in accordance with Conjecture \ref{conj31} below) that the
coarse moduli of $\Gg$-bundles has two connected components;
the neutral connected component is the classical Prym abelian variety 
${\rm ker}(1+\sigma^*)^0\subset {\rm Jac}(Y)$; here  $\sigma^*: {\rm Jac}(Y)\rightarrow {\rm Jac}(Y)$
is the induced involution on the Jacobian. 
Fix a line bundle $\M$ on $Y$ which satisfies ${\rm Norm}_{Y/X}(\M)\simeq \det(\pi_*(\O_Y))^{-1}$
(then $\sigma^*\M\simeq \M^{-1}$ and such a line bundle corresponds to $\M_0$ as above).
If $\L$ is a line bundle over $Y$ with ${\rm Norm}_{Y/X}(\L)\simeq \O_X$
(so that it corresponds to a $\Gg$-torsor),
then 
$$\begin{aligned}\det(\pi_*(\L\otimes \M))&\simeq \det(\pi_*(\O_Y))\otimes{\rm Norm}_{Y/X}(\L\otimes \M)\\
&\simeq \det(\pi_*(\O_Y))\otimes{\rm Norm}_{Y/X}(\L)\otimes{\rm Norm}_{Y/X}(\M)\simeq \O_X\ .
\end{aligned}$$
This shows that if $\L$ is a $\Gg$-torsor, the sheaf
 $\pi_*(\L\otimes\M)$ gives a ${\rm SL}_2$-bundle on $X$.

\subsection{}\label{pairings}
Suppose that ${\rm char}(k)\neq 2$ and that $\pi: \ti X\rightarrow X$ is a   (possibly ramified) double cover
with involution $\sigma$.
Consider the moduli stack of pairs of $(\E, \psi)$ of a ${\rm SL}_n$-bundle $\E$ over $\ti X\times_k S$
together with a perfect $\O_{X\times_kS}$-bilinear pairing
\begin{equation}
\psi: \pi_*(\E)\times \pi_*(\E)\rightarrow \pi_*(\O_{\ti X\times_k S})
\end{equation}
which is $\sigma$-hermitian in the sense that it satisfies 
$\psi(a\cdot v, w)=\psi(v, \sigma(a)\cdot w)$, $\psi(w,v)=\sigma(\psi(v, w))$
for $a\in\pi_*(\O_{\ti X\times_k S})$. Set $\Gg:={\rm SU}_n(\ti X/X)=({\rm Res}_{\ti X/X}{\rm SL}_n)^\sigma$;
here $\sigma$ acts on $g\in {\rm SL}_n(\O_{\ti X})$ by $g\mapsto J_n\cdot\sigma(g^{\rm tr})^{-1}\cdot J_n^{-1}$
where $J_n$ is the anti-diagonal unit matrix of size $n$.
Then the special unitary group $\Gg$ is a parahoric group scheme over $X$ and we can see that ${\mathcal M}_{\Gg/X}$
is the moduli stack of pairs above. This group scheme $\Gg$ is not of  the ``constant type" 
considered in \S \ref{constant}.

\section{}

We continue with the assumptions and notations of Section 2.
The first conjecture concerns the set of connected components, and is of Kottwitz style.  

\begin{conjecture}\label{conj31} Denote by $\pi_1(\mathcal G_{\bar\eta})$ the algebraic fundamental group of $\mathcal G_{\bar\eta}$ in the sense of Borovoi. Then
$$\pi_0(\mathcal M_{\mathcal G/X})\ =\ \pi_1(\mathcal G_{\bar\eta})_{\Gamma}.
$$
Here on the right hand side are the co-invariants under $\Gamma={\rm Gal}(\bar\eta/\eta)$. 
\end{conjecture}

\begin{Remarksnumb} {\rm In particular, if $\mathcal G_{\bar \eta}$ is semi-simple and simply connected, then $\mathcal M_{\mathcal G/X}$ should be connected. This would follow from Conjecture \ref{uniform} below and the fact that $LH$ is connected for any semi-simple simply connected group $H$ over $k((t))$, cf. 1.) in the Introduction. 
If $\mathcal G$ is constant, i.e comes by extension of scalars from a group scheme $G$ over $k$, then the action of $\Gamma$ on $\pi_1(\mathcal G_{\bar\eta})$ is trivial.  Over $\C$ the statement then follows from the topological uniformization theorem, \cite{S3}, Cor. 4.1.2. \qed
}\end{Remarksnumb}

The second conjecture concerns the uniformization of  $\mathcal M_{\mathcal G/X}$. 

\begin{conjecture}\label{uniform}
Let $x\in X(k)$. Let $\mathcal P$ be a $\mathcal G$-torsor over $X\times_k S$. 
If $\mathcal G_\eta$ is semi-simple, then after an $fppf$ base change $S'\rightarrow S$,
the restriction of $\mathcal P\times_S S'$ to $(X\smallsetminus \{x\})\times S'$ is trivial.
\end{conjecture}

Of course, one can also state a version of this conjecture involving a non-constant family of smooth connected projective curves, but this version would suffice to obtain a uniformization of  $\mathcal M_{\mathcal G/X}$. Namely, assuming $\mathcal G_\eta$ semi-simple, and choosing a uniformizer at $x$, we would have an isomorphism
\begin{equation}\label{quotient}
 \mathcal M_{\mathcal G/X}\ =\ \Gamma_{X\smallsetminus \{x\}}(\mathcal G)\backslash 
 L\mathcal G_{\eta_x}/L^+\mathcal G_x.
\end{equation}
Here $\Gamma_{X\smallsetminus \{x\}}(\mathcal G)$ denotes the ind-group scheme with $k$-rational points equal to 
$$
\Gamma_{X\smallsetminus \{x\}}(\mathcal G)(k)\ =\ 
\Gamma (X\smallsetminus \{x\}, \mathcal G)\ .
$$ More precisely, the expression (\ref{quotient}) represents the affine partial flag variety $\F_x= L\mathcal G_{\eta_x}/L^+\mathcal G_x$ as a $\Gamma_{X\smallsetminus \{x\}}(\mathcal G)$-torsor over  $ \mathcal M_{\mathcal G/X}$.  We will denote by $p_x$ the {\it uniformization morphism},
\begin{equation}\label{unif}
p_x: \F_x\to  \mathcal M_{\mathcal G/X}\ .
\end{equation}

\begin{Remarksnumb}{\rm  In the constant case $\mathcal G=G\times_{\Spec k}X$, this is the theorem of Drinfeld and Simpson [DS]. In the case $S=\Spec(k)$, the statement in this special case was proved much earlier by Harder [Ha]. Even in the context of the examples in  \ref{constant}, derived from constant group schemes, the conjecture is not trivial. \qed
}\end{Remarksnumb}

The third conjecture concerns the Picard group of $ \mathcal M_{\mathcal G/X}$. For this we assume that $\mathcal G_\eta$ is semi-simple, simply connected, and absolutely simple.   Let us also assume that $\mathcal G_{\eta_x}$ splits over a tamely ramified extension of $K_x=k((t))$. We recall from the theory of twisted loop groups [PR] that there is a natural homomorphism 
\begin{equation}
c_x:\  {\rm Pic}(\F_x)\rightarrow  \Z,
\end{equation}
the {\it central charge} (at $x$). Denoting by $X^*(\mathcal G(x))$ the character group of the fiber $\Gg(x)=\Gg\times_X\Spec(k(x))$ of $\mathcal G$ at $x$, we have an exact sequence
\begin{equation}\label{exact}
0\rightarrow X^*(\mathcal G(x))\rightarrow {\rm Pic}(\F_x) \xrightarrow {c_x} \Z\rightarrow 0\ ,
\end{equation}
which comes about as follows. There is a central extension $\tilde L\mathcal G_x$ of $L\mathcal G_x$ by $\mathbb G_m$ which acts on all line bundles on $\F_x$. Let $\tilde L^+\mathcal G_x$ be the restriction of this central extension to $L^+\mathcal G_x$. This defines a central extension $\tilde L^+\mathcal G_x$ of $L^+\mathcal G_x$ by $\mathbb G_m$, and an isomorphism
\begin{equation}\label{piciso}
X^*(\tilde  L^+\mathcal G_x )\xrightarrow{\ \sim\ } {\rm Pic}(\F_x) \ .
\end{equation}
 On the other hand, the reduction homomorphism $L^+\mathcal G_x\rightarrow \mathcal G(x)$  defines the exact sequence 
 \begin{equation}
 0\rightarrow X^*(\mathcal G(x))\rightarrow X^*(\tilde  L^+\mathcal G_x)\rightarrow  \Z\rightarrow 0\ ,
\end{equation}
which together with (\ref{piciso}) yields the exact sequence (\ref{exact}).   

Note that if $\mathcal G_x$ is  a special maximal parahoric group, then $\Gg(x)$ is an extension of a semi-simple group 
by a unipotent group, and so $X^*(\mathcal G(x))$ is trivial; this applies to all but finitely many points $x\in X(k)$. If  
$\mathcal G_x$ is a hyperspecial maximal parahoric group, then $\Gg(x)$ is semi-simple. Let us denote by
${\rm Bad} (\mathcal G)$ the set of points $x\in X(k)$ where $\mathcal G_x$ is not hyperspecial.

\begin{conjecture} Let $\mathcal G_\eta$ be semi-simple, simply connected and absolutely simple. We also assume that $\mathcal G_{\eta_x}$ splits over a tamely ramified extension of $K_x$, for all $x\in X(k)$. 

(i)  For any $x\in X(k)$, consider the homomorphism 
$$p_x^*: {\rm Pic} (\mathcal M_{\mathcal G/X})\rightarrow {\rm Pic}(\F_x)\ $$
induced by the uniformization morphism. 
Composing with $c_x$, we obtain a   homomorphism ${\rm Pic} (\mathcal M_{\mathcal G/X})\rightarrow\Z$.
If $x$ is not in ${\rm Bad} (\mathcal G)$ then this homomorphism is non-zero and independent of $x$.
Let us denote this homomorphism by $c$ or $c_{\mathcal G/X}$. 

(ii) Denote the  kernel of $c_{\mathcal G/X}$  by ${\rm Pic} (\mathcal M_{\mathcal G/X})^0$. There is a natural isomorphism
$${\rm Pic} (\mathcal M_{\mathcal G/X})^0\ \simeq\ \bigoplus_{x\in X(k)} X^*(\mathcal G(x)). $$

\end{conjecture}

\begin{Remarksnumb}\label{remarkHein}
{\rm  1) In the case that $\mathcal G=G\times_{\Spec k}X$, the 
 point (i) was proved by Sorger [S1]  for $k=\C$ and by Faltings [Fa1] for arbitrary $k$. In this case, (ii)  states that the homomorphism ${\rm Pic} (\mathcal M_{\mathcal G/X})\rightarrow {\rm Pic}(\F_x)$ is injective, which is also proved in these papers. 
In the case that $\mathcal G$ is derived from a constant group scheme, as described at the end of  section \ref{constant}, the point (ii) is proved by  Laszlo and Sorger in [LS]. 
 
2)  For any $x\in X(k)$, we can consider the homomorphism
$c_{\Gg/X, x}: {\rm Pic}(\M_{\Gg/X})\rightarrow \Z$ obtained by composing $p^*_x$ and $c_x$ as above.
As was pointed out by Heinloth (see also [He], Remark 15 (4)) we can have $c_{\Gg/X, x}\neq c_{\Gg/X}$
if  $x$ is in ${\rm Bad}(\Gg)$. Indeed, suppose that $\Gg={\rm SU}_{n}(\ti X/X)$ as in 
\S \ref{pairings} where $\pi: \ti X\to X$ is a {\sl ramified} double cover and suppose that $x\in X(k)$ is a 
branch point. Then $\Gg(\O_x)$ is a special but not hyperspecial parahoric subgroup. Set ${\mathcal H}={\rm Res}_{\ti X/X}{\rm SL}_{n}$. Then $\M_{{\mathcal H}/X}\simeq \M_{{\rm SL}_{n}/\ti X}$
and the Picard group of $\M_{{\mathcal H}/X}$ is isomorphic to $\Z$ with generator given by the determinant of
cohomology of the universal ${\rm SL}_{n}$-bundle over $\ti X$. Denote by $\delta$ the image of this element under 
 ${\rm Pic}(\M_{\mathcal H/X})\rightarrow{\rm Pic}(\M_{\Gg/X})$. Suppose that $y\in X(k)$ is {\sl not}
a branch point. Then $\mathcal H_y\simeq {\rm SL}_{n}\times {\rm SL}_{n}$, $\Gg_y\simeq {\rm SL}_{n}$,
with the embedding $\mathcal G_y\hookrightarrow \mathcal H_y$ given by $A\mapsto (A, J_n\cdot (A^{\rm tr})^{-1}\cdot J_n^{-1})$.
The corresponding map on Picard groups $ {\rm Pic}(\F_{{\mathcal H}_y}) \rightarrow {\rm Pic}(\F_{\Gg_y})$
is therefore given by the sum $\Z\times \Z\rightarrow \Z$. Hence, we can see that $c_{\Gg/X, y}(\delta)=2$. 
Assume now that $n$ is {\sl even}. Since $x$ is a branch point, by [PR] 10.4, the morphism $\F_{\Gg_x}\xrightarrow{} \F_{{\mathcal H}_x}$
given by $\Gg_x\hookrightarrow {\mathcal H}_x$ induces an isomorphism on Picard groups. Hence, we can see
that $c_{\Gg/X, x}(\delta)=1$.\qed
}\end{Remarksnumb}

We now come to the ``conformal blocks". Before this, we recall some facts from [PR] \S 10  about the Picard group of a partial affine flag variety $\F=LH/L^+P$. Here we are assuming that the group  $H$ over $k((t))$   
is semi-simple, simply connected and absolutely simple, and that $H$ splits over a tamely ramified extension of $k((t))$. Let $\{\alpha_i\mid i=1,\ldots,r\}$ be the set of affine roots corresponding to the walls bounding the facet in the Bruhat-Tits building fixed by the parahoric $P$. For each $i$ there is a closed embedding of a projective line into $\F$,
$$
\pione_{\alpha_i}\hookrightarrow \F \ .
$$
By associating to each line bundle on $\F$ the degree of its restriction to $\pione_{\alpha_i}$ for $i=1,\ldots,r$, we obtain an isomorphism (cf. [PR] Prop. 10.1),
\begin{equation}\label{deg}
{\rm deg}: {\rm Pic}(\F)\xrightarrow{\ \sim\ } \bigoplus\nolimits_{i=1}^r\Z\cdot \epsilon_i \ .
\end{equation}

 A line bundle $\L$ on $\F$ is called {\it dominant} if its image under (\ref{deg}) has all coefficients $\geq 0$.

Assume now that ${\rm char}(k)=0$. Then the Lie algebra of the universal extension $\tilde LH$ acts on the space of global sections ${\rm H}^0(\F, \L)$, and if $\L$ is dominant, this representation is the dual of the integrable highest weight representation corresponding to the element ${\rm deg}(\L)$. More precisely, we choose a minimal parahoric subgroup contained in $P$ (corresponding to an alcove in the Bruhat-Tits building containing the facet fixed by $P$) and set the  coefficients of all simple affine root $\alpha_i$ not occurring in (\ref{deg}) equal to zero. Then ${\rm deg}(\L)$ is a dominant integral weight $\lambda$ in the sense of Kac-Moody theory, and by Kumar and Mathieu \cite{Ku1},  \cite{Ku2}, \cite{Ma},
\begin{equation}
{\rm H}^0(\F, \L)=(V_\lambda)^* \ ,
\end{equation}
where on the RHS appears the dual of the integrable highest weight representation attached to $\lambda$. 

Let us spell out the above remarks in the standard case: assume that $H$ is constant, i.e.,  $H=G=G_0\times_kk((t))$. To simplify notations, let us assume also that the parahoric $P$ is an Iwahori subgroup which is contained in the maximal parahoric subgroup $P_0=G_0\times_kk[[t]]$. Hence $\F$ is the full affine flag variety. The target of the degree homomorphism can then be written in terms of  the fundamental weights,
\begin{equation}
{\rm deg}: {\rm Pic}(\F)\loniso \bigoplus\nolimits_{i=0}^l\Z\cdot\epsilon_i \ .
\end{equation}
Let us fix a maximal torus $T$ in the Borel subgroup of $G_0$ corresponding to $P$ in the sense of section \ref{constant}. We consider $T$ as a subgroup of $L^+P_0$, and let $\tilde T$ be the inverse image of $T$ in $\tilde LG$. Then the definition of $\tilde LG$ is such that each character of $\tilde T$ defines a line bundle on $\F$ and that in this way we obtain an isomorphism
\begin{equation}
X^*(\tilde T)\loniso{\rm Pic}(\F) \ .
\end{equation}
There is a unique splitting of the central extension $\tilde LG$ over $L^+P_0$. Hence we can  write canonically $\tilde T=T\times \G_m$, and 
\begin{equation}
X^*(\tilde T)=X^*(T)\oplus \Z \ .
\end{equation} 
In terms of this  decomposition the composed map  $X^*(\tilde T)\xrightarrow{\sim} {\rm Pic}(\F)\xrightarrow{\sim} \bigoplus \Z\cdot\epsilon_i$, is given as follows,
\begin{equation}\label{class}
\lambda=(\lambda^{(0)}, \ell) \mapsto \sum\nolimits_{i=1}^l n_i\cdot \epsilon_i+\Big(\ell-\sum\nolimits_{i=1}^ln_ir_i\Big)\cdot \epsilon_0 \ .
\end{equation}
Here $\lambda^{(0)}=\sum\nolimits_{i=1}^ln_i\epsilon_i$ and the positive integers $r_1,\ldots,r_l$ are the labels of the vertices of the dual  Dynkin diagram (denoted $a_i^\vee$ in Kac's book \cite{Kac}, p. 79). Note that $(\lambda^{(0)}, \ell)$ is dominant if and only if $\lambda^{(0)}$ is dominant and $\ell-\sum\nolimits_{i=1}^ln_ir_i\geq 0$. The last inequality can also be written in terms of the coroot $\theta^\vee$ for the highest root  in
 $X^*(T)_{\bf R}$. Indeed, $\theta^\vee=\sum\nolimits_{i=1}^l r_i\alpha_i^\vee$, in terms of the simple coroots $\alpha_1^\vee,\ldots,\alpha_l^\vee$. Hence the second condition for being dominant can be written in the familiar form, 
\begin{equation*}
\langle \theta^\vee, \lambda^{(0)}\rangle\leq \ell \ .
\end{equation*}

\medskip

We now return to the global situation and a general parahoric group scheme $\mathcal G$ over $X$. A   line bundle $\L$ on $\mathcal M_{\mathcal G /X}$ is called {\it dominant} if $p_x^*(\L)$ is a dominant line bundle on $\F_x$ for every $x$.  

\begin{conjecture}\label{conblock}
 Let ${\rm char}(k)=0$, and assume as before that $\Gg_\eta$ is semi-simple, simply connected and absolutely simple, and that $\Gg_{\eta_x}$ splits over a tamely ramified extension of $K_x$ for all $x\in X(k)$.  Let $S$ be a non-empty  finite subset of $X(k)$ containing 
${\rm Bad}(\mathcal G)$. Let $\L$ be a dominant line bundle on $\M_{\Gg/X}$
such that the central charge $c_{\Gg/X,x}(\L)$ is constant for $x$ in $X$.
There is a canonical isomorphism of finite-dimensional vector spaces
$$
{\rm H}^0(\mathcal M_{\mathcal G/X}, \L)\ \simeq \bigg [\bigotimes\nolimits_{x\in S}
{\rm H}^0(\F_x, p_x^*(\L))\bigg ]^{{\rm H}^0(X\setminus S,  {\rm Lie}( \mathcal G))} . 
$$
\end{conjecture}

Here the action of ${\rm H}^0(X\setminus S, {\rm Lie}( \mathcal G))$ comes from the fact that the homomorphism 
\begin{equation*}
{\rm H}^0(X\setminus S, {\rm Lie}( \mathcal G))\to \bigoplus\nolimits_{x\in S} {\rm Lie} (\mathcal G_{\eta_x})
\end{equation*}
lifts uniquely to the factor space of $\bigoplus\nolimits_{x\in S} {\rm Lie} (\tilde L \mathcal  G_{\eta_x})$ where the central elements in the central extensions for all $x\in S$ are  identified (here the assumption that $S\supset {\rm Bad} (\mathcal G)$ enters).  It is known \cite{B}, cor. 2.4, cf.\ also \cite{S1}, Prop.\ 2.3.2, 
that if $S$ is enlarged to $S'\supset S$, the RHS does not change.

\begin{Remarksnumb}{\rm
In the ``classical" theory, when $\Gg$ is constant, i.e., where $\mathcal G=G\times_{\Spec k}X$, one considers data which formally look very similar to the data above. Indeed, in the classical theory, just as here, one also fixes a  finite set $S$ of points, and  dominant integral weights, one for each point $x_i\in S$. These are written traditionally as in (\ref{class}) above in the form  $\lambda_i=(\lambda_i^{(0)}, \ell)$, where  $\lambda_i^{(0)}$ is a dominant weight for $G$ and $\ell\in\Z$ is the central charge with  $\langle \theta^\vee, \lambda_i^{(0)}\rangle\leq \ell$. These additional points and dominant integral weights are introduced to formulate and prove the fusion rules, which ultimately lead to an explicit determination of the dimension of the vector spaces in Conjecture  \ref{conblock}.

On the other hand, in [LS] the set $S$ and the dominant integral weights $\lambda_i$ appear for essentially the same reason as here (namely, to describe ${\rm H}^0(\F_{x_i}, p_{x_i}^*(\L))$ for $x_i\in {\rm Bad}(\Gg)$), except that here the situation is more general. In particular, in [LS], Thm. 1.2., the set $S$ consists of ${\rm Bad}(\Gg)$ and one additional point.  In \cite{B-L}, the parahoric group scheme $\Gg$ is the constant group scheme ${\rm SL}_n$ and the set ${\rm Bad}(\Gg)$ is empty, and $S$ consists of an arbitrary point of $X$, comp.\ also [B], Remarks in (2.6). 

Beauville in [B], Part I, treats formal properties of the spaces of conformal blocks which
appear in [LS] and only mentions in passing the geometric interpretation by the LHS in Conjecture \ref{conblock}. 
\qed

}\end{Remarksnumb}
In the classical case, when $\mathcal G=G\times_{\Spec k}X$,  the dimension of the RHS 
in Conjecture \ref{conblock} has been calculated by Faltings \cite{FalVerl}, \cite{FalCourse} by using the factorization rules and the fusion algebra, at least when $G$ is a classical group or of type $G_2$, comp.\ also [B], Part III.   It would be interesting to have a Verlinde type dimension formula in the  case of a general parahoric group scheme. Also, in the light of \cite{Te}, it should be possible to go beyond the case of dominant line bundles on $\mathcal M_{\mathcal G/X}$ and also consider higher cohomology groups.

\end{document}